%% file: acl_latex.tex
\pdfoutput=1
\documentclass[11pt]{article}

\usepackage{acl}
\usepackage{booktabs}
\usepackage{multicol}
\usepackage{multirow}
\usepackage{times}
\usepackage{latexsym}
\usepackage{amsmath}
\usepackage{amssymb}
\usepackage{xcolor}
\usepackage{pifont}
\usepackage{graphicx}
\usepackage{longtable}
\usepackage[T1]{fontenc}
\usepackage[utf8]{inputenc}

\usepackage{microtype}

\usepackage{inconsolata}

\title{From Good to Great: Improving Math Reasoning with Tool-Augmented Interleaf Prompting }

\author{
 Nuo Chen$^\clubsuit$
\quad
Hongguang Li$^{\diamondsuit}$  
\quad
{\bf Baoyuan Wang$^{\diamondsuit}$}  
{\bf \quad Jia Li$^\clubsuit$}\\
\\
  $^\clubsuit$Hong Kong University of Science and Technology (Guangzhou)\\ Hong Kong University of Science and Technology\\
  $^{\diamondsuit}$Xiaobing.AI\\
    \texttt{nchen022@connect.ust.hk}, \texttt{jialee@ust.hk}\\}

\begin{document}
\maketitle
\begin{abstract}
This paper investigates the performance of Large Language Models (LLMs) and Tool-augmented LLMs in tackling complex mathematical reasoning tasks. We introduce IMP-TIP: Improving Math Reasoning with Tool-augmented Interleaf Prompting, a framework that combines the strengths of both LLMs and Tool-augmented LLMs. IMP-TIP follows the ``From Good to Great" concept, collecting multiple potential solutions from both LLMs and their Tool-Augmented counterparts for the same math problem, and then selecting or re-generating the most accurate answer after cross-checking these solutions via tool-augmented interleaf prompting. The framework incorporates two key aspects: self-prompt and tool-augmented interleaf prompting (TIP). The former allows LLMs to autonomously refine and improve an initial prompt related to tool usage, while the latter enables LLMs to derive the final answer by dynamically analyzing the problem, cross-checking potential solutions, and revising previous reasoning hints in an interleaved manner. Experimental analysis shows that IMP-TIP achieves enhanced mathematical capabilities and outperforms traditional LLMs and tool-augmented LLMs in accuracy and reasoning diversity on math reasoning tasks. For instance, IMP-TIP can improve Tool-augmented ChatGPT on GSM8K-Hard from 56.0$\%$ to 65.2 $\%$.
\end{abstract}

\input{Sections/1Introduction}
\input{Sections/2relatedwork}
\input{Sections/3Method}
\input{Sections/4Experiments}
\input{Sections/5Analysis}
\input{Sections/6Conclusion}

% \input{sections/3Probing_Tasks}
% \input{sections/4Experiments}
% \input{sections/5Layer_wise_experiments}
% \input{sections/6XMATH_EXPERIMENTS}

% These instructions are for authors submitting papers to *ACL conferences using \LaTeX. They are not self-contained. All authors must follow the general instructions for *ACL proceedings,\footnote{\url{http://acl-org.github.io/ACLPUB/formatting.html}} and this document contains additional instructions for the \LaTeX{} style files.

% The templates include the \LaTeX{} source of this document (\texttt{acl.tex}),
% the \LaTeX{} style file used to format it (\texttt{acl.sty}),
% an ACL bibliography style (\texttt{acl\_natbib.bst}),
% an example bibliography (\texttt{custom.bib}),
% and the bibliography for the ACL Anthology (\texttt{anthology.bib}).

% \section{Engines}

% To produce a PDF file, pdf\LaTeX{} is strongly recommended (over original \LaTeX{} plus dvips+ps2pdf or dvipdf). Xe\LaTeX{} also produces PDF files, and is especially suitable for text in non-Latin scripts.

% Entries for the entire Anthology, followed by custom entries
\bibliography{anthology,custom}

\appendix

% \section{Example Appendix}
% 
\input{Sections/7Appendix}
% This is an appendix.

\end{document}

%% file: Sections/1Introduction.tex
\section{Introduction}

Large language models (LLMs) \cite{DBLP:journals/corr/abs-2005-14165, hu2021lora, zeng2022glm, openai2023gpt4,chen2023orca,you-etal-2022-end} such as ChatGPT have exhibited remarkable performances across a wide array of downstream tasks. They can effortlessly tackle downstream tasks by conditioning on a scant number of in-context exemplars or plain natural language task descriptions \cite{brown2020language}. Notwithstanding these significant advancements, even the most extensive LLMs are confronted with challenges when faced with intricate tasks that necessitate multiple reasoning steps \cite{gao2023pal, chen2023breaking}.

\input{Tables/pilot_analysis}

The capacity of LLMs to tackle complex tasks has been extensively assessed using math reasoning datasets such as GSM8K \cite{DBLP:journals/corr/abs-2110-14168} and SVAMP \cite{patel-etal-2021-nlp}. These datasets present math problems that cannot be directly answered but require multi-step reasoning. To encourage LLMs to engage in step-by-step reasoning, the chain-of-thought (COT) strategy \cite{wei2022chain} has emerged as the standard prompting approach. Through COT, LLMs are guided to generate a solution that consists of a sequence of intermediate steps,  leading to the final answer.
Nevertheless, it has been observed in previous studies that LLMs are prone to making mistakes or hallucinations, particularly during intermediate numerical computations \cite{qian2022limitations, yuan2023well, lu2022survey}. Even a minor mistake at this stage can result in a completely incorrect final answer. In an effort to address this limitation, a series of studies \cite{parisi2022talm, schick2023toolformer, yao2022react} have been undertaken, which leverage external tools such as the calculator to compensate for the weaknesses of LLMs. These approaches have shown significant improvements in the accuracy of answers on math reasoning tasks.

Although these approaches have led to significant improvements overall, a pertinent question arises: Can these tool-augmented LLMs outperform their traditional LLM-COT counterparts consistently? To explore this question,  we conduct a pilot analysis at first: we test both LLMs-COT and its tool-augmented counterpart (calculator-augmented) on GSM8K dataset, separately. Table \ref{pilot} meticulously compares the quantities of accurately and inaccurately predicted samples between the two methods. Upon careful analysis of the tabulated results, we can observe that incorporating a calculator has undeniably yielded a substantial boost (739 vs. 715) in the accuracy of  LLMs on the GSM8K dataset. This improvement can be primarily attributed to the calculator's capability to mitigate the potential errors that could arise during mathematical computations within the model. Nevertheless, it is noteworthy that despite this tool augmentation, there persist cases (\textbf{294}) where the LLMs-COT outperforms its Tool-augmented LLMs counterpart, achieving accurate predictions where the latter faltered. Upon experimental analysis, the primary factors contributing to this phenomenon can be categorized as:
1) LLMs may generate different reasoning logic for the same math question due to the varying token sampling probabilities. 2) When employing tool-augmented LLMs to tackle math reasoning problems, stringent requirements are imposed on the output text to adhere to specific formats. For instance, calculator-augmented LLMs \cite{schick2023toolformer} must produce the output text in a predefined numerical equation format, ensuring seamless compatibility with the calculator's invocation. Conversely, program-augmented LLMs \cite{gao2023pal} necessitate the generation of specific code functions that can be executed effectively. These prescribed format requirements may influence the mathematical reasoning process of LLMs.

\begin{figure*}
    \centering
    \includegraphics[width=0.96\linewidth]{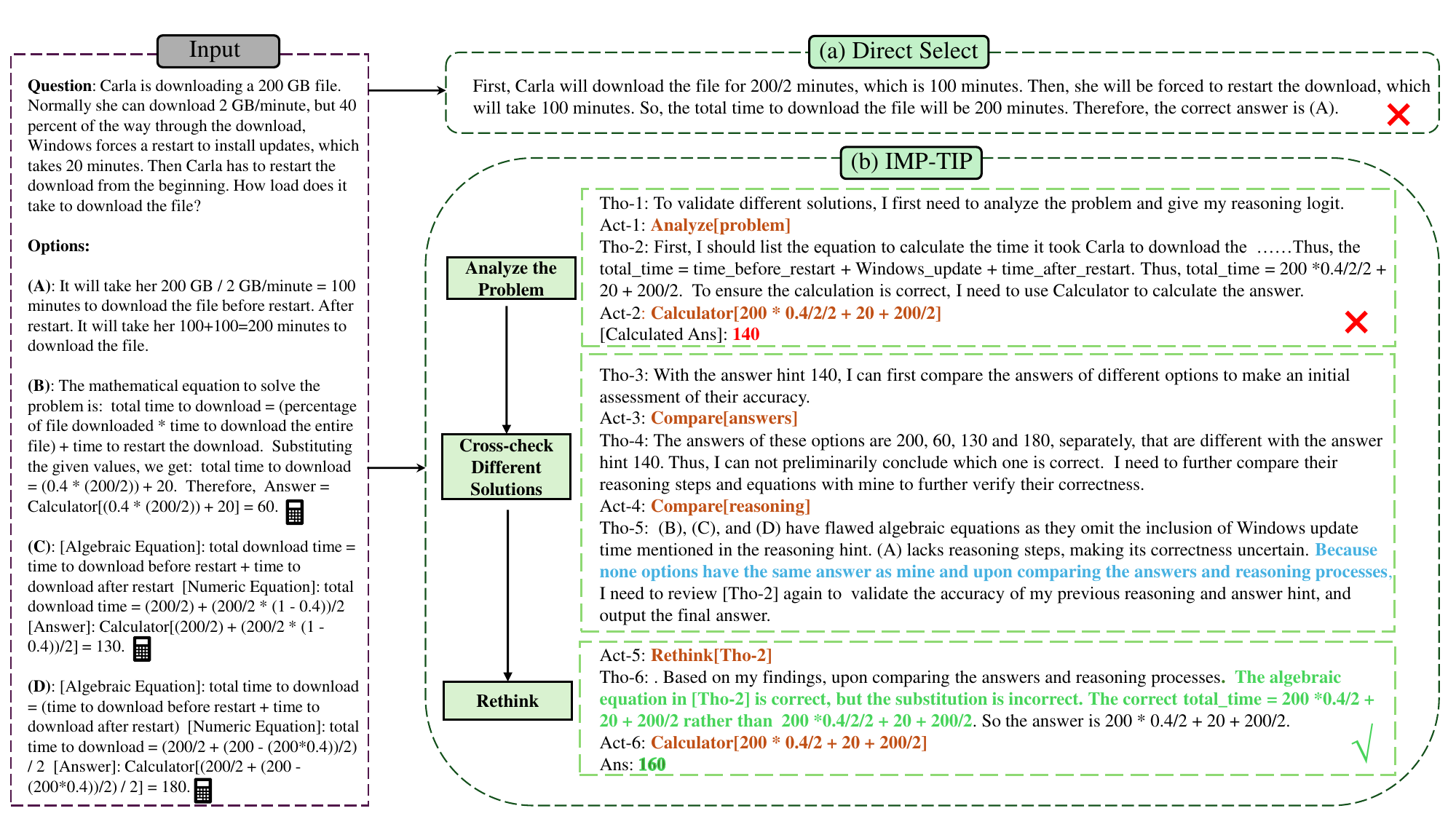}
    \caption{Comparison of 2 prompting methods: (a) \textit{Direct Select}, which requires LLMs to select the correct one from options.  (b) IMP-TIP, which derives the correct answer in an interleaved manner with \textit{action} and \textit{thought}. In this example, we present 3 options (Option  B, C, D) from tool-augmented LLMs (using the calculator) and 1 option from LLMs COT (Option A). Tho and Act are short for Thought and Action, separately.}
    \label{modeloverview}
    \vspace{-10pt}
\end{figure*}

Based on the aforementioned analysis, we aim to capitalize on the strengths of both LLMs and Tool-augmented LLMs to further enhance their mathematical capabilities. In this work, we introduce \textbf{IMP-TIP}: \textbf{I}mproving \textbf{M}ath \textbf{R}easoning with \textbf{T}ool-augmented \textbf{I}nterleaf \textbf{P}rompting. IMP-TIP is a framework that follows the idea ``\textit{From Good to Great}" that first collects multiple potential solutions from LLMs and their Tool-Augmented approaches given the same math problem, and then selects the most correct answer or re-generates new answer after cross-checking these solutions in a tool-augmented interleaved manner. To achieve this, we design IMP-TIP in two aspects: (1) We first propose \textit{self-prompt} to
address the challenge of crafting clear, diverse and effective tool-based prompts. It allows LLMs to autonomously refine and improve an initial prompt related to tool usage, resulting in enhanced prompts. Moreover, the self-prompt can yield multiple diverse tool-based prompts through iterative runs without requiring repetitive manual modifications,  eliciting varied reasoning paths (solutions) in LLMs. (2) Then given multiple solutions, we introduce \textit{tool-augmented interleaf prompting} (TIP), a paradigm combining reasoning, acting, and external tools in an interleaved manner to arrive at the correct answer. TIP allows LLMs to dynamically analyze the problem, cross-check potential solutions, and re-evaluate even revising previous reasoning hints through \textit{action} and \textit{thought}. Our technical contributions are summarized as follows:
\begin{itemize}
    \item We conduct an in-depth study for LLMs-COT and tool-augmented LLMs on the math reasoning tasks, make interesting observations and design a novel framework named IMP-TIP to improve LLMs.
    \item We introduce self-prompt to obtain diverse and machine-friendly tool-based prompts without manual effort.
    % without human-written effort. 
    \item We carry out extensive experiments to verify that our proposed IMP-TIP significantly improves LLMs and their tool-augmented approaches on 5 math reasoning tasks based on two different backbones.
\end{itemize}

% to increase the diversity of tool-based prompts in LLMs. Self-prompt allows LLMs to autonomously summarize and refine initial prompts related to tool usage, achieving improved and diverse tool-based prompts without the need for repetitive manual writing and modifications; (2) Through running self-prompt k times, we obtain k diverse tool-based prompts, resulting in k diverse solutions

%% file: Tables/pilot_analysis.tex
\begin{table}[]
\centering
\small
\begin{tabular}{c|c|c|c|c}

\toprule
Method&\multicolumn{4}{c}{Tool-Augmented  LLMs } \\
\midrule
\multirow{5}{*}{  LLMs-COT} & GSM8K &Wrong & Right &Total \\
\cmidrule{2-5}
% \midrule
&Wrong & \textit{286} & 318 & 604 \\
\cmidrule{2-5}
&Right&\textbf{294} & 421 &715\\
\cmidrule{2-5}
& Total &580&739&1319 \\
\bottomrule
\end{tabular}
\caption{Confusion matrix between LLMs and its Tool-Augmented approach on GSM8K dataset. Here, we select the calculator as the external tool, which is implemented by \texttt{eval} function in  \texttt{Python}. The experimental LLM is \texttt{text-davinci-003}. Prompts of two approaches contain two examples (2-shot) that are shown in Appendix \ref{prompts}.
}
\vspace{-15pt}
\label{pilot}

\end{table}

% \multirow{4}{*}{zero-shot} 

%% file: Sections/2relatedwork.tex
\section{Related Work}

\paragraph{Math Reasoning with LLMs} 
In recent times, large language models (LLMs) have demonstrated remarkable abilities in handling complex reasoning tasks \cite{scao2022bloom, DBLP:journals/corr/abs-2110-14168, zhou2022least, weng2022large, chen2023teaching,chen2023breaking,chen2022would}. Rather than providing direct final answers as outputs, prior research has shown that by employing diverse prompting methods such as Chain-of-Thought (COT) prompting \cite{wei2022chain}, LLMs can be guided through step-by-step reasoning, resulting in notably improved performance across a wide array of reasoning tasks. \citet{DBLP:conf/acl/ImaniD023} propose to generate multiple algebraic expressions or Python functions to solve the same math problem, aiming to explore different potential solutions. \citet{li2023making}
introduce a step-aware verifier to check the reasoning steps in COT, improving the reasoning capabilities. 
Self-Consistency \cite{wang2022self} is another effective work that combines different solutions and gets a final answer by aggregating to retrieve the most consistent answer. Among them, self-consistency bears resemblance to our proposed TIP, but our main distinction lies in the following: Given several solutions, our TIP can re-generate a new answer after cross-validating these solutions. In contrast, self-consistency can only select the most consistent answer from the existing ones.

% couples the generation of reasoning paths and a final answer by sampling from the decoder and
% using aggregation to retrieve the most consistent answer without extra modules.

\paragraph{Tool-Augmented LLMs} Currently, researchers have undertaken a wide array of studies aimed at enriching the step-by-step reasoning process.  These approaches include investigating the utilization of external tools \cite{parisi2022talm, schick2023toolformer, yao2022react}, like program interpreters \cite{lyu2023faithful, chen2022program} and the calculator \cite{schick2023toolformer}, training and utilizing external reasoning modules. For example, ReAct \cite{yao2022react} proposes an interleaved framework with utilizing an external search engine to solve multi-hop question answering tasks.
More recently, Toolformer \cite{schick2023toolformer} introduces a pipeline for training  LLMs that can call tools during training and inference. Parallel to these works, our work can be seen as a preliminary exploration of \textit{how to better utilize tools}, which involves a fusion of diverse solutions of LLMs and tool-augmented LLMs, culminating in the attainment of the final answer through tool-augmented interleaf prompting.

%% file: Sections/3Method.tex
\section{Methodogy}
In this section, we aim to illustrate our method in detail, as seen in Figure \ref{modeloverview} and Figure \ref{self_prompt}.
We first review the problem formulation
of math problem reasoning. Then we introduce our proposed self-prompt and tool-augmented interleaf prompting methods, sequentially.
\subsection{Problem formulation}
A math problem solving task can be defined as $\{Q, O, A\}$, where $Q$ is the target math question, $O = \{O_1, O_2,..., O_k \}$ are answer options if $Q$  is a K-way multiple choice problem, $A$ is the corresponding ground-truth answer. Given $Q$ and $O$ as inputs, LLMs can directly output answers or output a sequence of tokens as intermediate reasoning steps $R$ via COT. Then we can 
 obtain the answer in $R$ through regular expression matching.

% output a sequence of tokens or directly output

% \

\subsection{IMP-TIP}
The key insights of the proposed methods are two-fold: (1) We first propose self-prompt to obtain diverse tool-based prompts to induce more diverse reasoning paths and solutions. (2) We then introduce TIP to better derive the final answer from multiple reasoning solutions from LLMs and tool-augmented LLMs. The flow of our TIP can be defined as: \textbf{Initial answer determination from question analysis} $\mathbf{\rightarrow}$ \textbf{Cross-Validation of options accuracy} $\mathbf{\rightarrow}$ \textbf{Further verification of reasoning if discrepancies arise} $\mathbf{\rightarrow}$ \textbf{Rethinking for answer accuracy if uncertain} $\mathbf{\rightarrow}$ \textbf{Final answer}.

\subsubsection{Self-Prompt}
\begin{figure}
    \centering
    \includegraphics[width=0.95\linewidth]{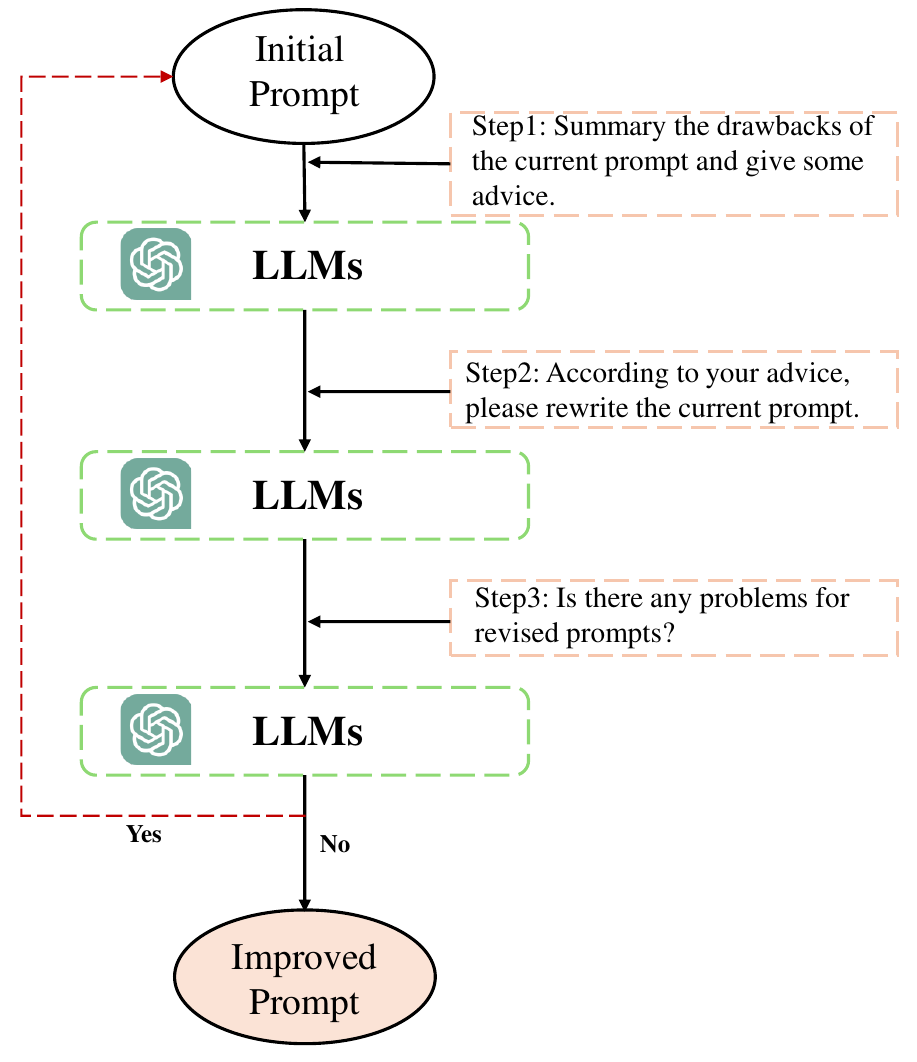}
    \caption{Overview of self-prompt.}
    \label{self_prompt}
    \vspace{-20pt}
\end{figure}

Crafting clear, diverse and concise tool-based prompts, which encompasses the output format, tool definition, and usage instructions, proves more challenging than traditional COT.  Thus, we propose self-prompt to address this issue. Moreover, different from previous works \cite{wang2022self,li2023making} that focus on generating k reasoning paths from one fixed prompt by sampling decoding, which limits the diversity of reasoning paths. Our proposed self-prompt can produce diverse prompts, automatically increasing prompt diversity and encouraging LLMs to think differently, eliminating the need for repetitive manual writing and editing.
% Previous works shows that 
% prove that exploring diverse reasoning paths can benefit the reasoning process. D O
Concretely, we design self-prompt in three steps, as shown in Figure \ref{self_prompt}:
\begin{itemize}
    \item Step1: We first give an initial tool-based prompt to the LLM with the instruction: ``Summary the drawbacks of the current prompt and give some advice''. 
    \item Step2: With the advice and prompt problems that are outputted from LLMs, we then instruct the LLM with: ``According to your advice, please rewrite the current prompt''. Thus, the LLM can output the revised prompt.
    \item Step3: Given the revised prompt as input and ask the LLM: ``Is there any problems for the revised prompt?''.  If the LLM responds ``Yes'' we will repeat steps 1 to 3 until the LLM answers ``No'' At that point, the revised prompt obtained in step 2 is the final improved prompt.
\end{itemize}

% The most simple implementation of self-prompt is to interact with ChatGPT or GPT-4 in OpenAI' website\footnote{chat.openai.com}.
Overall, our self-prompt follows the idea that ``LLMs know themselves better''.
Through the steps mentioned above, we not only address the issue of manually crafting reasonable and clear instructions for LLMs to use tools but also enhance the diversity of prompts simultaneously via running self-prompt through multiple iterations.  Experimentally, the improved tool-based prompt can help LLMs achieve better performances compared with the initial prompt.
We can run it M times, resulting in M diverse prompts, and then sample N reasoning paths for each prompt via sampling decoding. This way, we can obtain K = M $\times$ N tool-based diverse reasoning paths for each math problem.

% Following the line, , we introduce the self-prompt approach as a viable solution. The most obvious advantage of self-prompt is that it can be run multiple times, automatically increasing prompt diversity and encouraging LLMs to think differently, eliminating the need for repetitive manual writing and editing.

% we introduce 

% previous works focus on generating k reasoning paths from the fixed prompt by sampling decoding (1 to k).
% \citet{} use multiple human-written prompts to 

% by \

\subsubsection{Tool-augmented Interleaf Prompting} The core idea of TIP is: We augment the reasoning paths for each question from 1 to K+L, where  L from LLMs \footnote{In our experiments, we obtain L reasoning paths from LLMs by sampling decoding. } and K from their tool-augmented LLMs (obtained from self-prompt). Then we derive the final answer by observing and validating these solutions (reasoning paths) in an interleaved manner, which can be formalized as $ \mathcal{TIP}(\hat{A}|Q,R_1,...R_{K+L})$.
TIP is similar to how a human brain learns when solving mathematical and multiple-choice questions: given a question and different possible solutions, analyze the problem, observe the solutions to find how the question should be solved, and then learn from these solutions, memorize or revise its own solution, conclude the answer at last.

% think multi-
% ple times to derive different possible results, conclude on how the question should be solved, and
% then learn from or memorize its own solution.
% In particular, We imitate the human thinking process when by decomposing TIP into three major steps. 
\label{steps}
Following this, We decompose TIP into three major steps. 
Given a question and multiple potential solutions: (1) \textsc{Analyze the problem}: Before validating different answers, LLMs first need to analyze the problem and give their own reasoning logical via generating algebraic and numerical expression; (2) \textsc{Cross-check different solutions}: Then based on the above reasoning hint, LLMs can cross-check different solutions. Considering LLMs could occasionally conclude accurate answers when generating wrong reasoning paths, we break down this into two sub-steps: compare the numerical answers and then cross-verify intermediate steps, which explores whether each solution is both sound and accurate. (3) \textsc{Rethink}: Unlike traditional multiple-choice question answering tasks, where there is at least one correct solution and reliable reasoning hints, in our case, all four solutions could be incorrect, and previous reasoning hints might also be flawed. Hence, we introduce an extra step called ``Rethink'' which requires LLMs to review and correct their initial answer hints after observing and checking different solutions, outputting the final answer. 
We argue that LLMs are capable of gaining a deeper understanding of the problem by cross-checking various solutions.

Inspired by ReAct \cite{yao2022react}, we let LLMs achieve each step in an interleaved manner with \textit{thought} and \textit{action}. In each step, LLM receives the environmental context $c_t$ and takes an action $a_t \in \mathcal{A}$, where $c_t = \{t_1, a_1,...,t_{t-1},a_{t-1} \}$ and $\mathcal{A}$ is the action space. $t_1, a_1$ are previous \textit{thought} and \textit{action} at step 1. As shown in Figure \ref{modeloverview}, there could be multiple types of useful \textit{thoughts}, like creating action plans (\textit{Tho-1} in Figure \ref{modeloverview}), reasoning process of executing the action (\textit{Tho-2} in Figure \ref{modeloverview}), reasoning over the context (\textit{Tho-6} in Figure \ref{modeloverview}).

\paragraph{Action Space} As we focus on math reasoning tasks, $\mathcal{A}$ consists of the following actions to support interactive TIP: (1) \textit{Analyze[]}, which requires analyzing the specific solution or the math problem. (2) \textit{Compare[]}, which requires comparing and cross-checking answers or intermediate steps from different solutions or the generated answer hints. (3) \textit{Rethink[]}, which means to double-check or review the specific thinking process in previous \textit{thoughts}. (4) \textit{Calculator[]}, which calls the calculator to compute the expression.

Figure \ref{modeloverview} shows a typical example of IMP-TIP. In this example, all options are the wrong solutions, and the \textit{direct select} and \textit{self-consistency} approaches also give the wrong answers. 
In our proposed TMP-TIP, despite making an initial error in predicting during the ``analyze the problem" step, the model self-corrects by cross-checking and identifying the correct equation formulation in the \textit{Tho-2}, albeit with an algebraic mistake. In the \textit{Tho-6}, it rectifies the error, resulting in the correct answer via calling the calculator.

Overall, our IMP-TIP has the following salient features for solving math reasoning tasks: (I)
Multi-view Reasoning paths: Approaching the math problem from multiple perspectives from LLMs and tool-augmented LLMs; (II) Multi-verification: Comparing and cross-checking the different solutions to a reasoning hint, including verifying the correctness of numeric answers and the intermediate steps, aiming at providing a more accurate understanding of the problem; (III) Calculation-Verification: Using the calculator to calculate the expression or verify the calculation; (IV) Self-Check: Re-checking the reasoning paths that include algebraic equation and substitution in previous reasoning hints after cross-checking different solutions.

In theory, M, N, K and L could be assigned higher values, resulting in more reasoning paths. However, due to the fact that these K+L solutions will be subsequently used as inputs for our TIP, we must take into account the limitations posed by input token length and associated costs in the in-context learning setting. Experimentally, we have determined that setting M=3, N=1, K = M $\times$ N =3 and L=1 is appropriate.

% decompose TIP into three stages: (1)

%% file: Sections/4Experiments.tex
\section{Experiments}
\input{Tables/main_results}
\input{Tables/ablation}

In this section, we first review our testing datasets and evaluation metrics.
Then we briefly introduce our experimental settings that include backbones,  and  our prompts. At last, we will present our main results and ablation studies.

\subsection{Datasets and Metrics}
% In this paper, we validate our proposed approach on five math reasoning datasets: MAWPS, SVAMP, SVAMP-Hard, GSM8K, and GSM8K-Hard to show its effectiveness. 
% meticulously crafted by human problem writers to ensure high quality. 
\paragraph{Datasets} We validate our proposed approach on five math reasoning datasets:
(1) \textbf{MAWPS} \cite{koncel-kedziorski-etal-2016-mawps} serves as a collection of math word problems (MWPs), offering a unified testbed with 1921 samples for accessing various algorithms. (2) \textbf{SVAMP} \cite{patel-etal-2021-nlp}, short for Simple Variations on Arithmetic Math Word Problems, is an elementary-level MWP dataset that contains 1000 examples in test set. (3) \textbf{GSM8K} \cite{DBLP:journals/corr/abs-2110-14168} stands as a dataset comprising 1391 linguistically diverse grade school MWPs, meticulously crafted by human problem writers to ensure high quality. (4) \textbf{GSM8K-Hard} \cite{gao2023pal}, an advanced iteration of GSM8K, which substitutes numbers in the questions with larger counterparts. This modification aims to assess the generalization capability of LLMs when dealing with substantial numerical values. (5) Given the limited range of three-digit numbers in the SVAMP dataset, we enhanced model assessment by randomly replacing query numbers with values between 100,000 and 10,000,000. This adjustment, while maintaining original reasoning logic, lead to the creation of \textbf{SVAMP-Hard}. SVAMP-Hard and GSM8K-Hard provide a more rigorous test of LLMs' mathematical computational capabilities. 
% Table \ref{} shows the statistic of these datasets. 
% Notice that, we only use the test sets of them to validate the models' performances as we focus on in-context learning.

\paragraph{Metrics} Following previous standard works, we use the \textbf{Accuracy} as the evaluation metrics on all datasets.

\subsection{Experimental Settings}
\paragraph{Backbone} We conduct our experiments based on two OpenAI LLMs: GPT-3 (\texttt{text-davinci-003}) and ChatGPT (\texttt{gpt-3.5-turbo}) through Azure Service. We use the default parameters except temperature = 0 for greedy decoding during inference. In our experiments, we select the calculator as the external tool for tool-augmented approaches, which is implemented by  Python \textit{eval()} function.

\paragraph{Exemplars} All quintet mathematical reasoning datasets allocate two exemplars, which are arbitrarily selected solely from the GSM8K-trainset, for both LLMs-COT and tool-enhanced LLMs. In TIP, we manually compose trajectories to use three exemplars in the prompts. Each trajectory includes multiple thought-action steps, where each thought is free-form. LLMs combine these thoughts to devise a plan (``... I first need to analyze x''), guide the reasoning process (``Option (A) is right, I then ...''), call tool (``use the calculator to compute y''), and conclude the answer (``the answer is z''). The model will end with ``[Ans]'' when solving the task.
See Appendix A for more details about LLMs-COT, tool-augmented LLMs and TIP prompts.

\begin{figure*}
    \centering
    \includegraphics[width=0.9\linewidth]{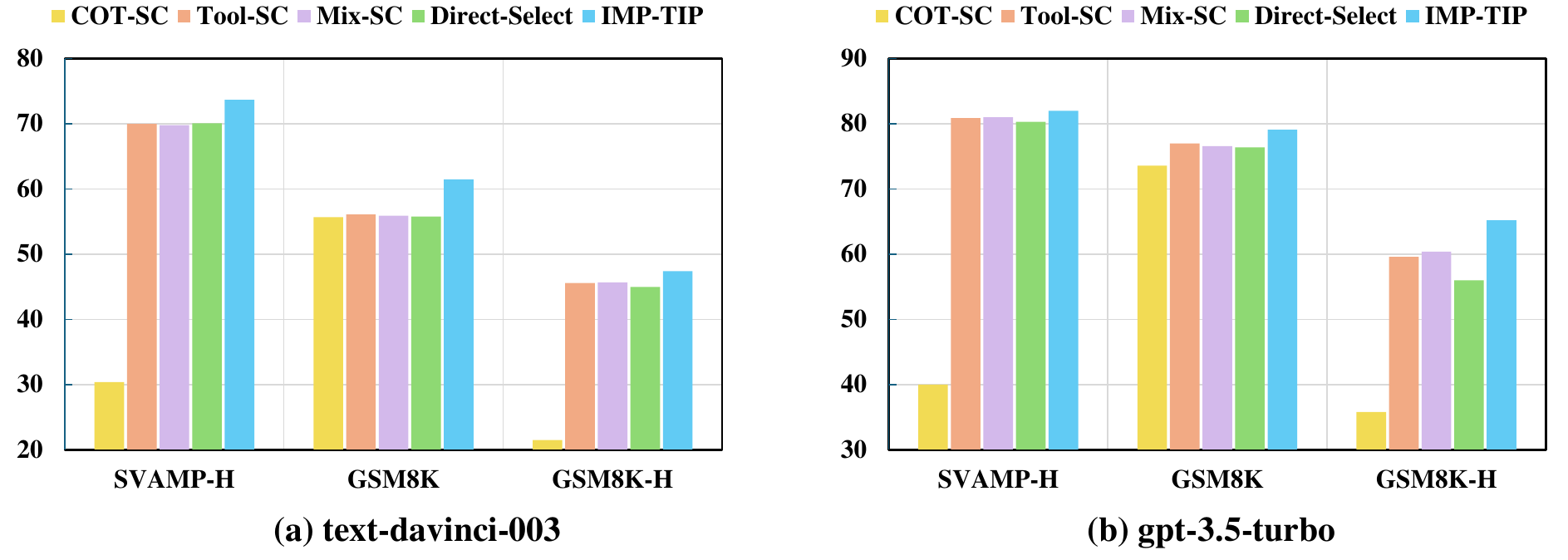}
    \caption{Ablation studies of different prompting methods based on \texttt{text-davinci-003} and \texttt{gpt-3.5-turbo}.}
    \label{consitency}
    \vspace{-15pt}
\end{figure*}

\paragraph{Baselines}
We mainly compare our approach with the following baselines: (1) \textbf{COT}, which solves the problem via step-by-step reasoning with our own implementations. (2) \textbf{Tool-Augmented}, signifies LLMs employing calculators for expression computation with \textit{initial tool-based prompts}, as elucidated in preceding sections. In this way, LLMs always need to follow a specific output format.
(4)
\textbf{Self-Prompt}, which utilizes the \textit{improved too-based prompts} from self-prompt to solve math reasoning tasks.
% also denominated as \textbf{tool-augmented LLMs}, signifies LLMs employing calculators for expression computation, as elucidated in preceding sections. In
% this way, LLMs always need to follow a specific output format.
% Tool-augmented LLMs which denote LLMs use the calculator to compute expression. In this way, LLMs always need to follow a specific output format.  In our experiments, tool-augmented prompts are obtained from self-prompt
(5) \textbf{Toolformer} \cite{schick2023toolformer}, a sophisticated paradigm that enables language models like GPT-J to utilize external tools, optimizing their performance in various tasks through fine-tuning. (4) \textbf{ART} \cite{paranjape2023art}, a  framework facilitating automatic multi-step reasoning and tool-use integration for LLMs.

For ablations, we also compare ours with (1) \textit{Self-consistency}: We build self-consistency baselines based on LLMs and their tool-augmented counterparts, named \textbf{COT-SC} and \textbf{Tool-SC}, separately. For a fair comparison, we sample reasoning paths 3 times by setting the temperature as 0.7. Then we also introduce \textbf{Mix-SC}, which votes the most consistent answer from the given solutions in TIP as the prediction.
(2) \textit{Direct Select}: Given a math problem, multiple solutions from LLMs and tool-augmented LLMs, we regard it as a multiple choice question answering task and require LLMs directly selecting one of them as the answer through COT. We show an example of it in Figure \ref{modeloverview}. See Appendix A for more details about their definitions.

\textbf{Of note}, although other works, such as PAL \cite{gao2023pal}, employs tools like programming languages for math problem-solving, we do not compare our method with theirs in this paper to maintain fairness, as our tools differ.

\subsection{Main Results}
Table \ref{tab_main_mlqa} shows the results of three backbone models with various prompting algorithms. Some key observations are as follows from the table:

\paragraph{LLMs still suffer in math calculations.} It is evident that the performances of these two backbone LLMs with COT on SVAMP-Hard and GSM8K-Hard are notably lower compared to their performance on SVAMP and GSM8K (e.g., 77.1$\%$ vs.29.5$\%$ on SVAMP and SVAMP-Hard). This discrepancy suggests that when tasked with computations involving larger numerical values, these LLMs exhibit an increased susceptibility to calculation errors. Interestingly, ChatGPT shows better computation ability than GPT-3.5.

\paragraph{Self-prompt benefits the performances for tool-using.} Across all datasets, it is discernible that the prompts subjected to self-prompt adjustments manifest enhanced performance relative to their original counterparts. This enhancement is attributed to the provision of more lucid and explicit tool usage instructions, alongside a heightened standardization of input-output formats.

\paragraph{IMP-TIP outperforms other baselines consistently, particularly on more challenging datasets involving complex computations.} Firstly, the table highlights the tool-augmented method's superiority over the traditional COT approach on all datasets, attributed to enhanced computational abilities. Secondly, our approach yields additional enhancements across all datasets, notably pronounced in the more challenging SVAMP-Hard and GSM8K-Hard sets demanding higher computational prowess. For example, IMP-TIP based on ChatGPT could improve COT and Self-Prompt from 35.0$\%$, 56.0$\%$ to 65.2$\%$ on GSM8K-Hard.

% \paragraph{}

\subsection{Ablation Study}

In this part, we conduct ablation studies to verify the impact of (1) each key component in IMP-TIP; (2) self-prompt (\textbf{How well prompts from self-prompt perform specifically?} See Appendix B for more details.).

\paragraph{Key Components}
% \subsubsection{Ablation Study}
% In this part, we conduct ablation studies to verify the impact of  each key component in IMP-TIP.
Given the sequential nature of our three-step IMP process (See Section 3), our ablation study comprises three sub-experiments in Table 3: 1) Step I only; 2) Step II only; 3) Steps I and II together. Table \ref{ablation} shows that: 1. Step II holds paramount importance; 2. Steps I and III also exhibit significance; 3. The combined use of all three steps yields optimal performance.

% 1. Step II holds paramount importance; 2. Steps I and III also exhibit significance; 3. The combined use of all three steps yields optimal performance.

% \paragraph{Impact of Self-Prompt}

% Then we employ 3 examples in the prompts for TIP.

%% file: Tables/main_results.tex
\begin{table*}[t]
\small
    \begin{center}
    \setlength{\tabcolsep}{5pt}
    \begin{tabular}{llcccccc}
    
    \toprule
  \multirow{2}{*}{\textbf{Backbone}}&  \multirow{2}{*}{\textbf{Algorithms}} & \multicolumn{5}{c}{\textbf{Dataset}} \\ \cmidrule{3-7}

  &&\textbf{MAWPS} & \textbf{SVAMP} & \textbf{SVAMP-H} & \textbf{GSM8K} & \textbf{GSM8K-H} &  \textbf{Average}\\
    \midrule
\texttt{GPT-J}&Toolformer & 44.0 & 29.4 & - &- &-&-\\

\midrule
    \multirow{5}{*}{\texttt{\parbox{0.2\linewidth}{
GPT-3\newline(text-davinci-003)
}}} &ART & 90.1 & 76.2 & -&-&-&-\\
    & COT & 89.6 
    & 77.1 & 29.5 & 54.2 & 20.4 & 54.2 \\
    & Tool-Augmented & 89.9 & 77.6 & 69.9 & 55.5 & 44.1 &  67.3\\
   & Self-Prompt & 90.1 & 78.5 & 70.0 & 56.1 & 44.5 &  67.8\\
   & \textbf{IMP-TIP} & \textbf{90.7} & \textbf{79.3} & \textbf{73.7} & \textbf{61.5} & \textbf{47.4} &  \textbf{70.5}\\

    \midrule
    \multirow{4}{*}{\texttt{\parbox{0.2\linewidth}{
ChatGPT\newline(gpt-3.5-turbo)
}}} & COT & 91.0 & 76.5 & 39.2 & 73.0 & 35.0 & 62.9 \\
& Tool-Augmented & 91.9 & 77.0 & 78.5 & 75.5 & 54.0 &  75.4\\
   & Self-Prompt & 92.1 & 78.0 & 79.7 & 76.0 & 56.0 &  76.4\\
   & \textbf{IMP-TIP} & \textbf{92.6} & \textbf{79.3} & \textbf{82.0} & \textbf{79.1} & \textbf{65.2} &  \textbf{79.6}\\
%  &  \textbf{Ours-LRSM}  & 80.38\,/\,66.46 & \textbf{70.34}\,/\,\textbf{52.35} & \textbf{65.84}\,/\,\textbf{50.34} & \textbf{59.58}\,/\,\textbf{42.11} & \textbf{64.07}\,/\,\textbf{47.88}& \textbf{69.26}\,/\,\textbf{48.79} &  \textbf{68.25}\,/\,\textbf{51.32} \\
%  &  \quad w/ \textit{Single-Language} & 80.20\,/\,66.68 & 70.10\,/\,51.40 & 64.99\,/\,49.81 & 58.19\,/\,40.50 & 62.40\,/\,46.58 & 68.11\,/\,47.47 &  67.43\,/\,50.41 \\
%  &  \quad w/ \textit{Cross-lingual}(en) & \textbf{80.40}\,/\,66.81 & 70.23\,/\,51.79 & 65.21\,/\,50.01 & 59.04\,/\,41.52 & 63.66\,/\,47.15 & 69.09\,/\,48.41 &  67.95\,/\,50.95 \\

    \bottomrule
    \end{tabular}
    
    \end{center}
\caption{The overall  results on the five datasets.  We highlight the best results for each backbone. In our experiments, as we obtain 3 improved tool-based prompts through self-prompt (Section 3), we report their average results. SVAMP-H and GSM8K-H are short for SVAMP-Hard and GSM8K-Hard. We report average results of 3 runs.}\label{tab_main_mlqa}
\vspace{-15pt}
\end{table*}

%% file: Tables/ablation.tex
\begin{table}[]
\small
\centering

\begin{tabular}{l|ccc}

\toprule

\textbf{Models}  & SVAMP-H &GSM8K& GSM8K-H \\ \midrule

\textbf{IMP-TIP} & \textbf{73.7} & \textbf{61.5} & \textbf{47.4} \\ 
\midrule
w/ Step-I &70.4 & 56.3 & 44.5 \\

w/ Step-II &72.4 & 58.7 & 45.6 \\
\midrule
w/ Step-I + II &72.9 & 59.4 & 46.5 \\

% w/o Step-I&72.4 & 58.7 & 45.6 \\
% w/o Step-II
% w/o Step-III&72.9 & 59.4 & 46.5 \\
% \midrule
% \multicolumn{2}{l}{\textbf{Re-generate new answers}} &3.1$\%$& 8.1$\%$\\
\bottomrule
\end{tabular}
\caption{Ablation studies in our IMP-TIP. Here, we conduct results based on \texttt{text-davinci-003}.}
\vspace{-15pt}
\label{ablation}

\end{table}

%% file: Sections/5Analysis.tex
\section{Analysis}
\input{Tables/Case}
\paragraph{IMP-TIP vs. Self-consistency} Given multiple solutions from LLMs and tool-augmented LLMs, two intuitive ways to derive the answer are: 
% 1) \textit{Self-Consistency}, selecting the most consistent answer among solutions as the final answer; 2) \textit{Direct Select}, regarding it as a multiple question answering task and directly selecting one of them as the answer.
1) \textit{Self-Consistency}; 2) \textit{Direct Select}.
As aforementioned, we have three versions of self-consistency: COT-SC, Tool-SC and Mix-SC. Figure \ref{consitency} shows the comparison between ours and these approaches on three datasets. We can observe that Tool-SC and Mix-SC perform comparably but they perform much better than COT-SC, just as tool-augmented LLMs surpass LLM-COT. Across two distinct backbone models, IMP-TIP consistently outperforms other methods on each dataset.

\input{Tables/propotion}

\subsection{Why IMP-TIP works?}
In pursuit of a comprehensive dissection of our IMP-TIP's underlying mechanics, we proffer a distributional analysis capturing the trajectory through which IMP-TIP reaches its definitive answers. Table \ref{proportion} delineates the proportions in which IMP-TIP either selects answers from provided options or re-generates new answers on SVAMP-Hard and GSM8K-Hard datasets. The former is segmented into three categories of assessment by IMP-TIP: solely based on LLM's answer, based on both the LLM and the Tool-augmented methods, and exclusively within the Tool-augmented method. Predominantly, IMP-TIP favors answers from Tool-augmented LLMs that are derived from provided solutions. Notably, IMP-TIP's capacity for generating new answers through observation and cross-checking significantly enhances its performance relative to self-consistency. This unique capability distinguishes IMP-TIP from other methods.

\input{Tables/rebuttal_analysis}

We remain curious about the performance of IMP-TIP when not all the given options are necessarily correct. To verify this, we conduct in-depth analysis of the GSM8K-Hard dataset in Table \ref{rebuttal_1} shows that our model incorrectly answers 2.7\% of cases when all provided answer choices are correct. Conversely, it correctly answers 5.2\% of cases when all choices are incorrect. \textbf{We have provided corresponding cases in Appendix, Table 5-6.} These outcomes highlight the challenges in enabling a model to self-verify and self-correct. Despite some errors, \textbf{our method strengthens LLM's ability to enhance its understanding and correct mistakes through interleaving reasoning steps}, which is a central aspect of our approach in the in-context learning setting.

\paragraph{Case Study} Table \ref{tab:case1} shows the an example from IMP-TIP and self-consistency on the GSM8K dataset. 
In this instance, only option (D) stands correct, with option (A) exhibiting a computational error, and both options (B) and (C) featuring accurate algebraic expressions albeit flawed numerical substitutions. Within this context, Mix-SC would favor the most consistently occurring answer, leading to the selection of 44 as the predicted answer, despite its inaccuracy. While our approach also encounters numerical substitution errors in [Tho-2], meticulous scrutiny of both the answer and reasoning steps reveals the preceding mistake, affirming the correctness of option (D). This case underscores another inherent advantage of IMP-TIP over self-consistency: the capacity to derive an accurate answer through thoughtful reasoning, rather than a simplistic reliance on majority voting for consistency. We present more cases in Appendix C.

% and error analyses
% \paragraph{Error Analysis}  

%% file: Tables/Case.tex
\begin{table*}[!t]\footnotesize
\centering
\small

\begin{tabular}{p{0.95\linewidth}}
% \begin{tabular}{lc}
% \multirow{28}{*}{\textbf{Persona Prompts}} & 
\toprule
\textbf{Question}: Sadie slept 8 hours on Monday. For next two days, she slept 2 hours less, each, because she had to complete some assignments. If the rest of the week she slept 1 hour more than those two days, how many hours did she sleep in total throughout the week? \\
% \textbf{Options}: \\
\textbf{(A)}: Sadie slept 8 hours on Monday. For the next two days, she slept 2 hours less, each, so she slept 6 hours each day. If the rest of the week she slept 1 hour more than those two days, she slept 7 hours each day. So the total amount of hours she slept throughout the week is 8+6+6+7+7+7+7=42. So the answer is 42. \\
\textbf{(B)}:  The mathematical equation to solve the problem is:  total hours slept = (Monday's hours + (2 days * (hours - 2)) + (4 days * (hours + 1))).  Substituting the given values, we get:  total hours slept = (8 + (2 * (6 - 2)) + (4 * (6 + 1))).  Therefore,  Answer = Calculator[8 + (2 * (6 - 2)) + (4 * (6 + 1))] = 44. \\
\textbf{(C)}:  [Algebraic Equation]: Total hours slept = (Monday's hours + (2 days * (Monday's hours - 2 hours)) + (4 days * (Monday's hours - 2 hours + 1 hour))  [Numeric Equation]: Total hours slept = (8 + (2 * (8 - 2)) + (4 * (8 - 2 )))  [Answer]: Calculator[(8 + (2 * (8 - 2)) + (4 * (8 - 2 )))] = 44. \\
\textbf{(D)}:  [Algebraic Equation]: total hours slept = hours slept on Monday + (hours slept on Tuesday + hours slept on Wednesday + hours slept on Thursday + hours slept on Friday + hours slept on Saturday + hours slept on Sunday)  [Numeric Equation]: total hours slept = 8 + (6 + 6 + 7 + 7 + 7 + 7)  [Answer]: Calculator[8 + (6 + 6 + 7 + 7 + 7 + 7)] = 48. \\
\midrule
\textbf{IMP-TIP}: \\
 Tho-1: To validate different solutions, I first need to analyze the math problem and give my reasoning logic.\\
Act-1: Analyze[problem] \\
 Tho-2: I need to find the total amount of hours Sadie slept throughout the week given the information that she slept 8 hours on Monday, 2 hours less for the next two days, and 1 hour more for the rest of the week. Thus, the total hours = 8 + (2 * (8 - 2)) + (4 * (8 - 2)). To ensure the calculation is correct, I need to use Calculator to calculate the answer. \\
Act-2: Calculator[8 + (2 * (8 - 2)) + (4 * (8 - 2 ))]\\
 Tho-3: With the answer hint 44, I can  first compare their answers to make an initial assessment of their accuracy.
 \\
 Act-3: Compare[answers] \\
 Tho-4: I observe that  (A) directly gives the answer while (B), (C) and (D) use the Calculator to compute the answer. The answers of them are 42, 44, 44, 48. Half of them align with the answer hint 44, I can preliminarily conclude that (B), (C) are correct. But there are still two different answers, I cannot determine which answer is definitely correct. Thus, I can compare their reasoning steps and equations with mine to further verify their correctness.
 \\
 Act-4: Compare[reasoning] \\
 Tho-5: (A) and (D) follow similar reasoning paths and their mathematical equations are the same, but their answers are different. Due to (D) using the calculator, (A) must have mistakes in calculation.  (B) and (C) give similar algebraic equations, but their substitutions are incorrect. (B) states that Sadie slept 6 hours each day for the next two days, which is not accurate. (C) is also incorrect. I may also make a mistake in [Tho-2] because my numerical equation is the same as (C). I need to review [Tho-2] again to validate the accuracy of my previous reasoning and answer hint. \\
 Act-5: Rethink[Tho-2]. \\
 Tho-6: Based on my findings, the algebraic equation in [Tho-2] is correct, but the substitution is incorrect.  The correct total hours = 8 + (2 * (8 - 2)) + (4 * (8 - 2 + 1)) rather than 8 + (2 * (8 - 2)) + (4 * (8 - 2)).\\
Act-6: Calculator[8 + (2 * (8 - 2)) + (4 * (8 - 2 + 1))] \\
Tho-7: After revising the substitution in [Tho-2], I think the final answer is 48 and (D) is correct. \\
Ans: \textbf{48} \\
% \bottomrule
\midrule
\textbf{Mix-SC}: 44. \\
\bottomrule
\caption{Examples of generated solutions by IMP-TIP and self-consistency on GSM8K dataset.  } 
\label{tab:case1}
\end{tabular}

\vspace{-10mm}
\end{table*}

%% file: Tables/propotion.tex
\begin{table}[]
\centering
\small
\begin{tabular}{l|lcc}

\toprule

&\textbf{Answer Type}  & SVAMP-H & GSM8K-H \\ \midrule
\multirow{3}{*}{\textbf{Select}}&Only LLM& 4.5$\%$ & 6.3$\%$ \\ 
&LLM and Tool & 13.4$\%$ & 18.5$\%$ \\ 
&Only Tool & 79.0$\%$ & 67.1$\%$ \\
\midrule
\multicolumn{2}{l}{\textbf{Re-generate new answers}} &3.1$\%$& 8.1$\%$\\
\bottomrule
\end{tabular}
\caption{Proportion of each answer type in our IMP-TIP. Here, we conduct results based on \texttt{gpt-3.5-turbo}.}
\vspace{-20pt}
\label{proportion}

\end{table}

%% file: Tables/rebuttal_analysis.tex
\begin{table}[]
\centering
\small
\begin{tabular}{c|c|c|c|c}

\toprule
Method&\multicolumn{4}{c}{Answer options have the right answer? } \\
\midrule
\multirow{5}{*}{IMP-TIP} & GSM8K-H (\%) &Yes & No &Total \\
\cmidrule{2-5}
% \midrule
&Wrong & \textbf{2.7\% }& 31.8\% & 34.5\% \\
\cmidrule{2-5}
&Right&60.3\% &\textbf{5.2\%}  &65.5\%\\
\cmidrule{2-5}
& Total &63.0\%&37.0\%&100\% \\
\bottomrule
\end{tabular}
\caption{ Our
IMP-TIP's performance when provided with answer choices that are either correct or incorrect in GSM8K-Hard dataset. The experimental LLMs is \texttt{gpt-3.5-turbo}.}
\vspace{-20pt}
\label{rebuttal_1}

\end{table}

% \multirow{4}{*}{zero-shot} 

%% file: Sections/6Conclusion.tex
\section{Conclusion}

% In this paper, we first present a pilot analysis comparing the performance of traditional and tool-augmented LLMs on math reasoning tasks. The results show that while tool-augmentation improves accuracy, there are instances where traditional LLMs outperform their tool-augmented counterparts. 
% Motivated by this, 
This work introduces IMP-TIP, a framework designed to capitalize on the strengths of both conventional and tool-augmented LLMs, enhancing their mathematical capabilities. IMP-TIP employs self-prompting to refine its initial prompts. Additionally, it utilizes tool-augmented interleaf prompting (TIP) to reach the correct answer. This process involves interleaved \textit{Action} and \textit{Thought}, which encompasses analyzing the problem, observing, cross-verifying different solutions, and even rethinking previous answer hints. Through rigorous experimentation and comprehensive analysis of math reasoning tasks, IMP-TIP demonstrates a substantial enhancement over established baselines.

%% file: Sections/7Appendix.tex
\section{Appendix A: Prompts}
\label{prompts}
\input{Appendix/cot_prompt}
\begin{figure}
    \centering
    \includegraphics[width=0.95\linewidth]{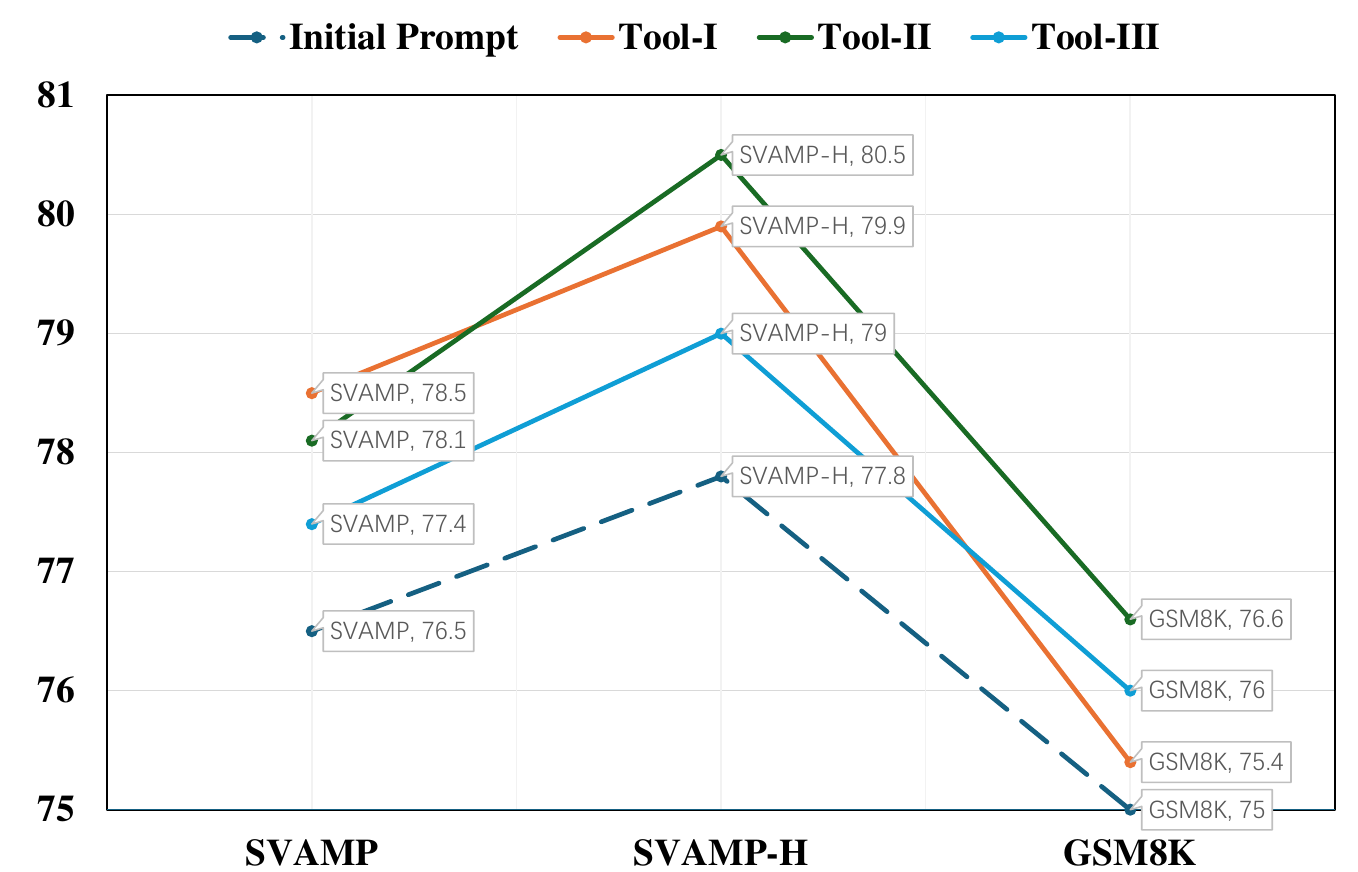}
    \caption{Model performances of different prompts. Here, we use the \texttt{gpt-3.5-turbo} as the backbone LLM.}
    \label{tools}
    \vspace{-15pt}
\end{figure}

\input{Appendix/tool_prompts}
Table \ref{tab:cot} presents COT prompts in our experiments.
Table \ref{tab:tool} shows three tool-based prompts obtained from our proposed self-prompt. 
It is worth noting that in most cases, the self-prompt does not alter the examples within the prompt. Its primary focus lies in modifying instructions related to tool usage, input-output formats, and tool definitions. All of them contain only 2 examples.

\input{Appendix/imp_prompt}

Table \ref{tab:imp} shows the prompts of TIP in our experiments. Practically, we append ``------'' as the stopping symbols after each \textit{Action} to achieve interleaving.

\subsection{Definition of Self-consistency and Direct-Select}

\begin{itemize}
    \item COT-SC: Utilizing COT prompts in Table 1 for math reasoning and sampling 3 times.
    \item COT-SC: Utilizing tool-based prompts in Table 2 for math reasoning and sampling 3 times. Then we report their average results.
    \item Mix-SC: Selecting the most consistent answer from given options (A, B, C, D) in TIP.
    \item Directly Select: Select the answer from the given options in TIP through COT. Prompts are shown in Table 1.
\end{itemize}

\section{Appendix B: Self-Prompt}
In this component, we conduct ablation studies to answer the question: How well do prompts from self-prompt perform? Figure \ref{tools} depicts a comparison between the initial tool-based prompt and the three prompts modified through self-prompt. It is evident across all three datasets that the prompts refined by self-prompt exhibit improved performance compared to the original initial prompts. This enhancement can be attributed to the provision of clearer, more explicit tool usage instructions and more standardized input-output formats. Overall, our self-prompt follows the road ``LLMs know themselves better.'' The most simple implementation of self-prompt is to interact with ChatGPT or GPT-4 in OpenAI' website\footnote{chat.openai.com}.
\section{Appendix C: More Cases}

\input{Appendix/easy_example}
\input{Appendix/negative}
\input{Appendix/false_positive}
In this section, we present typical cases when IMP-TIP solves math reasoning tasks. 

Table \ref{tab:easy_exm} presents the most straightforward scenario, wherein the provided solution corresponds with the answer generated by IMP-TIP in [Tho-2]. In this instance, the LLM detects the uniformity of all five answers upon comparison, resulting in a notable elevation of confidence and the subsequent direct issuance of the final answer.

Table \ref{tab:negative} also provides an example wherein IMP-TIP makes an error. Instead of further comparing the reasoning steps of different solutions after analyzing the answers, the model directly presents a conclusion, asserting the correctness of its approach.

 Table \ref{tab:false_pos} introduces another intriguing case. In this instance, IMP-TIP attains the correct answer for [Tho-2]. However, upon reevaluating the preceding reasoning hints and formulas for [Tho-6], it identifies a numerical substitution error within [Tho-2]. Despite this recognition, IMP-TIP proceeds to output the identical formula as the final answer. Instances in which LLMs exhibit self-contradictory reasoning yet still yield the correct ultimate answer in IMP-TIP, as exemplified here, are referred to as ``FALSE POSITIVE."

%% file: Appendix/cot_prompt.tex
\begin{table*}[!t]\footnotesize
\centering
\small

\begin{tabular}{p{0.95\linewidth}}
% \begin{tabular}{lc}
% \multirow{28}{*}{\textbf{Persona Prompts}} & 
\toprule
COT prompts \\
\midrule
Your task is to answer the following math questions.
    \\
Question: julia played tag with 18 kids on monday . she played tag with 10 kids on Tuesday. how many more kids did she play with on monday than on Tuesday? \\
Answer: Let's think step by step. julia playsed tag with 18 kids on monday and 10 kids tuesday, separately. So the amount of kids that she played with on monday than on tuesday is 18-10=8. So the answer is 8.\\
\\
Question: Jack had 9 action figures and 10 books on a shelf in his room. later he added 7 more action figures to the shelf. how many more action figures than books were on his shelf?\\
Answer: Let's think step by step. The amount of action figures that Jack had is 9+7=16. And Jack had 10 books. So the amount of action figures than books on his shelf is 16-10=6. So the answer is 7.\\
\midrule
\midrule
Initial tool-based prompts \\
\midrule
Your task is to solve the following middle-school arithmetic problems by using the Calculator.
A calculator is a system used for performing mathematical calculations, ranging from basic arithmetic to more complex functions. \\
    
You can do so by writing a mathematical equation and generating the answer format
starting with "Answer = Calculator[expression]", where "expression" is the expression to be computed. \\
\\
Below are examples: \\
Question: julia played tag with 18 kids on Monday. she played tag with 10 kids on Tuesday. how many more kids did she play with on monday than on Tuesday? \\
Thought:  The mathematical equation to solve the problem is:
number of additional kids played with = number of kids played with on Monday - number of kids played with on Tuesday.
Substituting the given values, we get:
number of additional kids played with = 18 - 10.
Therefore,  Answer = Calculator[18 - 10] \\
\\
Question: Jack had 9 action figures and 10 books on a shelf in his room. later he added 7 more action figures to the shelf. how many more action figures than books were on his shelf?\\
Thought: The mathematical equation to solve the problem is:
difference in number of action figures and books = (number of action figures initially + number of action figures added) - number of books.
Substituting the given values, we get:
difference in number of action figures and books = (9 + 7) - 10.
Therefore,  Answer = Calculator[( 9 + 7 )- 10] \\
\midrule
\midrule
Direct-select prompts \\
\midrule
There is a multiple choice question answering task. You are required to select one of the options to answer the problem. Let's think step by step. \\
% \midrule
\bottomrule
    \caption{LLMs-COT, initial tool-based prompts and direct-select prompts in our experiments.  } 
\label{tab:cot}
\end{tabular}

\vspace{-10mm}
\end{table*}

%% file: Appendix/tool_prompts.tex
\begin{table*}[!t]\footnotesize
\centering
\small

\begin{tabular}{p{0.95\linewidth}}
% \begin{tabular}{lc}
% \multirow{28}{*}{\textbf{Persona Prompts}} & 
\toprule
Prompt of Tool-I for math reasoning tasks. \\
\midrule

Your task is to solve the following middle-school arithmetic problems using a calculator. Follow the steps provided to write out the algebraic and numeric equations, and use a calculator to find the answers. You can generate the answer format starting with ``Answer = Calculator[expression]", where ``expression" is the expression to be computed. \\

Example Problems:\\

Question:
Julia played tag with 18 kids on Monday and 10 kids on Tuesday. How many more kids did she play with on Monday than on Tuesday?\\

Thought:\\

Algebraic Equation: Additional kids played with = Kids played with on Monday - Kids played with on Tuesday. \\
Numeric Equation: Additional kids played with = 18 - 10\\
Answer: Calculate Calculator[18 - 10]\\
\\
Question :
Jack had 9 action figures and 10 books on a shelf. Later, he added 7 more action figures to the shelf. How many more action figures than books were on his shelf?
\\
Thought:\\

Algebraic Equation: Difference in action figures and books = (Initial action figures + Added action figures) - Number of books.\\
Numeric Equation: Difference in action figures and books = (9 + 7) - 10\\
Answer: Calculate Calculator[(9 + 7) - 10]\\

Remember to enter the appropriate equation into your calculator to find the answer. You can also explain how you would verify your solution. \\

\midrule 
\midrule
Prompt of Tool-II for math reasoning tasks. \\
\midrule
Your task is to solve the following middle-school arithmetic problems by using the Calculator.
A calculator is a system used for performing mathematical calculations, ranging from basic arithmetic to more complex functions. \\
\\
You can do so by writing out Algebraic Equation, Numeric Equation, Answer steps. 
In the Answer step, you should generate the output format  with ``Calculator[expression]", where 
 ``expression" is the expression to be computed.\\

Below are examples: \\

Question: julia played tag with 18 kids on monday. she played tag with 10 kids on Tuesday. how many more kids did she play with on monday than on Tuesday? \\
Thought:
[Algebraic Equation]: number of additional kids played with = number of kids played with on Monday - number of kids played with on Tuesday. \\

[Numeric Equation]: number of additional kids played with = 18 - 10. \\

[Answer]: Calculator[18 - 10] \\

\\

Question: Jack had 9 action figures and 10 books on a shelf in his room. later he added 7 more action figures to the shelf. how many more action figures than books were on his shelf? \\

Thought:
[Algebraic Equation]: difference in number of action figures and books = (number of action figures initially + number of action figures added) - number of books.\\

[Numeric Equation]: difference in number of action figures and books = (9 + 7) - 10\\

[Answer]: Calculator[( 9 + 7 )- 10]\\
\midrule
\midrule
Prompt of Tool-III for math reasoning tasks. \\
\midrule

Your task is to solve the following middle-school arithmetic problems by using the Calculator.
A Calculator is a system used for performing mathematical calculations, ranging from basic arithmetic to more complex functions.  \\
\\

You can do so by writing out Algebraic Equation, Numeric Equation, and Answer steps:\\
(1) [Algebraic Equation], which directly builds an Algebraic Equation by using variables from the question to generate the solution.\\
(2) [Numeric Equation], which plugs in values from the question to transform the above Algebraic Equation into its numeric form.\\
(3) [Answer], which uses Calculator to calculate the above numeric equation, but not directly to give the final answer. You can generate the answer format with ``Calculator[expression]", where ``expression" is the expression to be computed.\\
\\

Below are examples: \\

Question: julia played tag with 18 kids on monday. she played tag with 10 kids on tuesday. how many more kids did she play with on monday than on Tuesday? \\

Thought:
[Algebraic Equation]: number of additional kids played with = number of kids played with on Monday - number of kids played with on Tuesday. \\
 
[Numeric Equation]: number of additional kids played with = 18 - 10. \\

[Answer]: Calculator[18 - 10] \\

\\
Question: Jack had 9 action figures and 10 books on a shelf in his room. later he added 7 more action figures to the shelf. how many more action figures than books were on his shelf? \\

Thought:
[Algebraic Equation]: difference in the number of figures and books = (number of figures initially + number of figures added) - number of books. \\

[Numeric Equation]: difference in the number of figures and books = (9 + 7) - 10\\

[Answer]: Calculator[( 9 + 7 )- 10]\\
You should keep in mind: \\
    (1) you can not image or output the result of "Calculator[expression]". \\
    (2) When generating a mathematical formula, do not include any units of measurement. \\

\bottomrule
\caption{Three tool-based prompts from Self-Prompt.  } 
\label{tab:tool}
\end{tabular}

\vspace{-10mm}
\end{table*}

%% file: Appendix/imp_prompt.tex
\begin{table*}[!t]\footnotesize
\centering
\tiny

\begin{tabular}{p{0.95\linewidth}}
% \begin{tabular}{lc}
% \multirow{28}{*}{\textbf{Persona Prompts}} & 
\toprule
Your task is to validate the several options to the same math problem and select the correct ones.  You should make your own reflections, revise your own thoughts, give the final answers based on the previous findings at the end. You can use the Calculator to help you calculate the numeric equations and verify the solutions. \\
A Calculator is a system used for performing mathematical calculations, ranging from basic arithmetic to more complex functions. \\
To do this, you can writing out interleaving Thought, Action and Answer steps. Notice that, you need to end with Answer when solving the task. \\
Below are examples:

Question: A robe takes 2287720 bolts of blue fiber and half that much white fiber.  How many bolts in total does it take?\\
(A): A robe ...... 2287720 + (2287720/2) = 3,431,280 bolts.\\
(B): Based on the statement of the question ...... total number of bolts = Calculator[2287720 + (2287720/2)] = 3431580. \\
(C): [Algebraic Equation]: total number of bolts = number of blue bolts + number of white bolts...... Calculator[2287720 + (2287720 / 2)] = 3431580. \\
(D): [Algebraic Equation]: total number of bolts ......Calculator[2287720 + (2287720/2)] = 3431580.  \\
    
Thoughts: 
[Thought-1]: To validate different solutions, I first need to analyze the math problem and give my reasoning logic. \\

[Action-1]: Analyze[problem] \\

[Thought-2]: First, the number of bolts of while fiber is half of the number of bolts of blue fiber. Total number of bolts = Number of bolts of blue fiber + Number of bolts of white fiber. Thus, total number of bolts = 2287720 + 2287720/2. To ensure the calculation is correct, I should use Caculator to calculate the answer. \\

[Action-2]: Calculator[2287720 + 2287720/2] \\

[Calculated Result]: 3431580 \\

[Thought-3]:  With the answer hint 3431580, I can first  compare the answers of different options to make an initial assessment of their accuracy. \\

[Action-3]: Compare[answers] \\

[Thought-4]: The answer of option (A) is 3431280 that contradicts the answer hint 3431580, while options (B), (C), and (D) align with the answer hint. Therefore, I can preliminarily conclude that (A) is incorrect, and (B), (C), and (D) are correct. Since most of the answer options are consistent with my previous answer hint, I believe the current answer hint is likely correct.  I can compare their reasoning steps and equations with mine to further verify my judgment and investigate why option (A) yielded an incorrect answer. \\

[Action-4]: Compare[reasoning] \\

[Thought-5]: Option (A) follows a similar reasoning process and formula as mine in [Thought-2], suggesting an error in the final calculation. Options (B), (C), and (D) provide equations that align with my earlier reasoning. Based on my findings, upon comparing the answers and reasoning processes, options (B), (C), and (D), along with the answer hint, are all correct. \\

[Answer]: Calculator[2287720 + 2287720/2] \\
    \\
    Question: Trent is 5 years older than Jane, and Jane is 3 years younger than Quinn. If Quinn is 2538570, how old is Trent? \\
    (A): Quinn is 2538570...... So the answer is Trent is 2538572. \\
    (B): The mathematical equation to solve the problem is:  ......Answer = Calculator[2538570 + (3 - 2538570) + 5] = 8. \\
    (C): [Algebraic Equation]: Trent's age = Quinn's age ......2538578. \\
    (D): [Algebraic Equation]: ...... Calculator[2538570 + 5 + 3] = 2538578. \\
    
    Thoughts:
    [Thought-1]: To validate different solutions, I first need to analyze the math problem  and give my reasoning logic. \\
    
    [Action-1]: Analyze[problem] \\
    % ------\\
    
    [Thought-2]: I need to find Trent's age given the relationships between Trent, Jane, and Quinn's ages. Jane is 3 years younger than Quinn and 5 year younger than Trent.  If Quinn is 2538570, Trent's age is 2538570 -3 + 5. To ensure the calculation is correct, I need to use Caculator to calculate the answer. \\
    
    [Action-2]: Calculator[2538570 -3 + 5] \\

    [Calculated Result]: 2538572 \\
    
    [Thought-3]: With the answer hint 2538572, I can frist compare the answers of different options to make an initial assessment of their accuracy. \\
    
    [Action-3]: Compare[answers] \\

    [Thought-4]: The answer of option (A) is 2538572 that align with the answer hint 2538572. However options (B)'s answer is 8, the answers of (C) and (D) are both 2538578 which all contradict the answer hint. Therefore, I can preliminarily conclude that (A) is correct, and (B), (C), and (D) are incorrect. But  there are three different answers and only one aligns with the answer hint, I cannot determine which answer is definitely correct. Thus, I can compare their reasoning steps and equations with mine to further verify their correctness. \\
    
    [Action-4]: Compare[reasoning] \\

    [Thought-5]: (A) follows a similar reasoning process and formula as mine in [Thought-2]. (B) is correct in terms of the mathematical equation used, but the substitution is incorrect. It states that Jane's age is equal to 3. (C) and (D) present that Trent's age = Quinn's age + 8 (5+3), meaning that Trent is 5 years older that Jane, and Jane is 3 year older than Quinn which contradict the information given in the question. The correct algebraic equation would be to subtract 3 from Quinn's age to determine Jane's age and then add 5 to Jane's age to find Trent's age. Based on my findings, I think (B), (C) and (D) are wrong, the answer hint and (A) are correct with no need improvements. I don't need to review my previous answer hint.\\
    
    [Answer]: Calculator[2538570 -3 + 5] \\
    \\
    Question: Carla is downloading a 200 GB file. Normally she can download 2 GB/minute, but 40 percent of the way through the download, Windows forces a restart to install updates, which takes 20 minutes. Then Carla has to restart the download from the beginning. How load does it take to download the file? \\
    (A): It will take her 200 minutes to download the file. \\
    (B): The mathematical equation ...... Therefore,  Answer = Calculator[(0.4 * (200/2)) + 20] = 60. \\
    (C): [Algebraic Equation]: ...... Calculator[(200/2) + (200/2 * (1 - 0.4))/2] = 130. \\
    (D): [Algebraic Equation]: total time ......[Answer]: Calculator[(200/2 + (200 - (200*0.4))/2) / 2] = 180. \\
    
    Thoughts:
    [Thought-1]: To validate different solutions, I first need to analyze the math problem and give my reasoning logic. \\
    
    [Action-1]: Analyze[problem] \\

    [Thought-2]: First, I should list the equation to calculate the time it took Carla to download the data before the Windows restart: time before restart = Downloaded data/Download speed. And time after restart = File size/Donwload speed. Thus, the total time = time before restart + Windows update + time after restart. Thus, total time = 200 *0.4/2/2 + 20 + 200/2.  To ensure the calculation is correct, I need to use Calculator to calculate the answer. \\
    
    [Action-2]: Calculator[200 * 0.4/2/2 + 20 + 200/2] \\
    % ------\\
    
    [Calculated Result]: 140 \\
    
    [Thought-3]: With the answer hint 140, I can first compare the answers of different options to make an initial assessment of their accuracy. \\
    
    [Action-3]: Compare[answers] \\
    
    [Thought-4]: The answers of these options are 200, 60, 130 and 180, separately, that are different with the answer hint 140. Thus, I can not preliminarily conclude which one is correct.  I need to further compare their reasoning steps and equations with mine to further verify their correctness.\\
    
    [Action-4]: Compare[reasoning]\\
    % -------\\
    
    [Thought-5]: (B), (C), and (D) have flawed algebraic equations as they omit the inclusion of Windows update time mentioned in the reasoning hint. (A) lacks reasoning steps, making its correctness uncertain. Because none options have the same answer as mine, I need to review [Thought-2] again to  validate the accuracy of my previous reasoning and answer hint.\\
    
    [Action-5]: Rethink[Thought-2]\\
    % ------\\
    
    [Thought-6]: The algebraic equation in [Thought-2] is correct, but the substitution is incorrect. The correct total time = 200 *0.4/2 + 20 + 200/2 rather than  200 *0.4/2/2 + 20 + 200/2. So based on my findings, the answer is Calculator[200 * 0.4/2 + 20 + 200/2].\\
    
    [Answer]: Calculator[200 * 0.4/2 + 20 + 200/2]\\
    \\
    Keep in mind the following requirements: \\
    (1) You must output the correct equation as the final answer in [Answer]. \\
    (2) If answers of all options are the same as the answer hint, you can directly end the task and output the correct equation in Answer. \\
    (3) If you decide to use Calculator in Action, generate the format with "Calculator[expression]", where "expression" is the expression to be computed. \\
    (4) Do not include any units of measurement in "Calculator[expression]". \\
    (5) You can not image the results of "Calculator[expression]". \\
\bottomrule
\caption{Prompts of TIP. In this table, we omit details of each option. } 
\label{tab:imp}
\end{tabular}

\vspace{-10mm}
\end{table*}

%% file: Appendix/easy_example.tex
\begin{table*}[!t]\footnotesize
\centering
\small

\begin{tabular}{p{0.95\linewidth}}
% \begin{tabular}{lc}
% \multirow{28}{*}{\textbf{Persona Prompts}} & 
\toprule
\textbf{Question}: There were some birds sitting on the fence. 725067 more birds came to join them. if there are a total of 544650 birds on the fence now how many birds had been sitting on the fence at the start? \\
% \textbf{Options}: \\
\textbf{(A)}:  The total number of birds on the fence now is 544650. And 725067 more birds came to join them. So the total number of birds at the start is 544650 - 725067 = -180417.  \\
\textbf{(B)}:  Answer = Calculator[544650 - 725067] = -180417 \\
\textbf{(C)}:  [Algebraic Equation]: number of birds at the start = total number of birds on the fence now - number of birds that came to join them.  [Numeric Equation]: number of birds at the start = 544650 - 725067  [Answer]: Calculator[544650 - 725067] = -180417 \\
\textbf{(D)}:  [Algebraic Equation]: number of birds at the start = total number of birds on the fence now - number of birds that came to join them.  [Numeric Equation]: number of birds at the start = 544650 - 725067  [Answer]: Calculator[544650 - 725067] = -180417
\\
\midrule
\textbf{IMP-TIP}: \\
 Tho-1: To validate different solutions, I first need to analyze the math problem and give my reasoning logic.\\
Act-1: Analyze[problem] \\
 Tho-2: The problem states that there were some birds sitting on the fence, and 725067 more birds came to join them. The total number of birds on the fence now is 544650. To find the number of birds that had been sitting on the fence at the start, we can subtract the number of birds that came to join from the total number of birds now. To ensure the calculation is correct, I need to use the Calculator to calculate the answer.  \\
 
Act-2:  Calculator[544650 - 725067] \\

 Tho-3: With the answer hint -180417, I can first compare the answers of different options to make an initial assessment of their accuracy. \\
 Act-3: Compare[answers] \\
 Tho-4: The answers of options (A), (B), (C), and (D) are all -180417, which align with the answer hint -180417. Therefore, I can conclude that all options are correct. No more action are needed.
 \\

Ans: \textbf{-180417.0} \\
% \bottomrule
\midrule
\textbf{Ground-Truth Anaswer}: \textbf{-180417.0} \\
\bottomrule
\caption{The \textsc{most simple} interleaved example of generated solutions by IMP-TIP on SVAMP-Hard dataset.  } 
\label{tab:easy_exm}
\end{tabular}

\vspace{-10mm}
\end{table*}

%% file: Appendix/negative.tex
\begin{table*}[!t]\footnotesize
\centering
\small

\begin{tabular}{p{0.95\linewidth}}
% \begin{tabular}{lc}
% \multirow{28}{*}{\textbf{Persona Prompts}} & 
\toprule
\textbf{Question}: julia played tag with 829557 kids on monday and 853729 kids on tuesday. she played cards wtih 913524 kids on Wednesday. how many kids did she play tag with altogether?\\
% \textbf{Options}: \\
\textbf{(A)}:  Julia played tag with 829,557 kids on Monday and 853,729 kids on Tuesday. To find out how many kids she played tag with altogether, we need to add these two numbers together.   829,557 + 853,729 = 1,683,286  So, Julia played tag with a total of 1,683,286 kids altogether. \\
\textbf{(B)}:  Answer = Calculator[829557 + 853729] = 1683286 \\
\textbf{(C)}:  [Algebraic Equation]: total number of kids played tag with = number of kids played tag with on Monday + number of kids played tag with on Tuesday + number of kids played tag with on Wednesday.  [Numeric Equation]: total number of kids played tag with = 829557 + 853729 + 913524  [Answer]: Calculator[829557 + 853729 + 913524] = 2596810 \\
\textbf{(D)}:  [Algebraic Equation]: total number of kids played tag with = number of kids played tag with on Monday + number of kids played tag with on Tuesday + number of kids played tag with on Wednesday.  [Numeric Equation]: total number of kids played tag with = 829557 + 853729 + 913524  [Answer]: Calculator[829557 + 853729 + 913524] = 2596810
\\
\midrule
\textbf{IMP-TIP}: \\
 Tho-1: To validate different solutions, I first need to analyze the math problem and give my reasoning logic.\\
Act-1: Analyze[problem] \\
 Tho-2: Julia played tag with a certain number of kids on Monday, Tuesday, and Wednesday. To find out how many kids she played tag with altogether, I need to add the number of kids she played tag with on each day. The equation to calculate the total number of kids played tag with is: total number of kids played tag with = number of kids played tag with on Monday + number of kids played tag with on Tuesday + number of kids played tag with on Wednesday. To ensure the calculation is correct, I need to use the Calculator to calculate the answer.\\
 
Act-2:  Calculator[829557 + 853729 + 913524] \\

 Tho-3: With the answer hint 2596810, I can first compare the answers of different options to make an initial assessment of their accuracy. \\
 Act-3: Compare[answers] \\
 Tho-4: The answers of options (A), (B), (C), and (D) are 1,683,286, 1,683,286, 2,596,810, and 2,596,810, respectively. The answer hint is 2,596,810. Therefore, options (A) and (B) are incorrect, while options (C) and (D) are correct. Since options (C) and (D) have the same answer and reasoning process, I can conclude that the correct answer is Calculator[829557 + 853729 + 913524].
 \\
 % Tho-5: Based on my findings, the correct answer is Calculator[829557 + 853729 + 913524] = 2596810
 % \\

Ans: \textbf{2596810} \\
% \bottomrule
\midrule
\textbf{Ground-Truth Anaswer}: \textbf{1683286.0} \\
\bottomrule
\caption{The typical \textsc{Negative} interleaved example of generated solutions by IMP-TIP on SVAMP-Hard dataset.  } 
\label{tab:negative}
\end{tabular}

\vspace{-10mm}
\end{table*}

%% file: Appendix/false_positive.tex
\begin{table*}[!t]\footnotesize
\centering
\small

\begin{tabular}{p{0.95\linewidth}}
% \begin{tabular}{lc}
% \multirow{28}{*}{\textbf{Persona Prompts}} & 
\toprule
\textbf{Question}: Sadie slept 8 hours on Monday. For next two days, she slept 2 hours less, each, because she had to complete some assignments. If the rest of the week she slept 1 hour more than those two days, how many hours did she sleep in total throughout the week? \\
% \textbf{Options}: \\
\textbf{(A)}: Wendi feeds each of her chickens three cups of mixed chicken feed, containing seeds, mealworms and vegetables to help keep them healthy. So the amount of cups of feed that Wendi needs to give her chickens in the final meal of the day is 20x3=60. So the answer is 60. \\
\textbf{(B)}:  The mathematical equation to solve the problem is:  total cups of feed = (number of chickens * cups of feed per meal) - (cups of feed given in morning + cups of feed given in afternoon).  Substituting the given values, we get:  total cups of feed = (20 * 3) - (6887483 + 25).  Therefore,  Answer = Calculator[(20 * 3) - (6887483 + 25)] = -6887448. \\
\textbf{(C)}:  [Algebraic Equation]: number of cups of feed for the final meal = (number of chickens * cups of feed per chicken) - (cups of feed given in the morning + cups of feed given in the afternoon).  [Numeric Equation]: number of cups of feed for the final meal = (20 * 3) - (6887483 + 25)  [Answer]: Calculator[(20 * 3) - (6887483 + 25)] =-6887448 \\
\textbf{(D)}:  [Algebraic Equation]: total cups of feed for the final meal = (number of chickens * cups of feed per chicken) - (cups of feed given in the morning + cups of feed given in the afternoon)  [Numeric Equation]: total cups of feed for the final meal = (20 * 3) - (6887483 + 25)  [Answer]: Calculator[(20 * 3) - (6887483 + 25)] = -6887448
\\
\midrule
\textbf{IMP-TIP}: \\
 Tho-1: To validate different solutions, I first need to analyze the math problem and give my reasoning logic.\\
Act-1: Analyze[problem] \\
 Tho-2: Wendi needs to give her chickens three cups of feed in three separate meals. She has already given 6887483 cups of feed in the morning and 25 cups of feed in the afternoon. Thus, the number of cups of feed that Wendi needs to give her chickens in the final meal of the day is (Number of chickens * Cups of feed per chicken) - (Cups of feed given in the morning + Cups of feed given in the afternoon). To ensure the calculation is correct, I need to use Calculator to calculate the answer. \\
 
Act-2: \textbf{Calculator[(20 * 3) - (6887483 + 25)]}\\

 Tho-3: With the answer hint -6887448, I can first compare the answers of different options to make an initial assessment of their accuracy. \\
 Act-3: Compare[answers] \\
 Tho-4: \textbf{The answers of these options are 60, -6887448, -6887448 and -6887448, separately, which are different with the answer hint -6887448.} Thus, I can not preliminarily conclude which one is correct. I need to further compare their reasoning steps and equations with mine to further verify their correctness.
 \\
 Act-4: Compare[reasoning] \\
 Tho-5: Option (A) follows a similar reasoning process as mine in [Tho-2], but the equation is incorrect. Options (B), (C), and (D) provide equations that align with my earlier reasoning. However, the substitution is incorrect as it does not include the cups of feed given in the morning and afternoon. Based on my findings, upon comparing the answers and reasoning processes, options (A), (B), (C), and (D) are all incorrect. I need to review [Tho-2] again to validate the accuracy of my previous reasoning and answer hint. \\
 Act-5: Rethink[Tho-2]. \\
 Tho-6: The algebraic equation in [Tho-2] is correct, but the substitution is incorrect. The correct total cups of feed = (Number of chickens * Cups of feed per chicken) - (Cups of feed given in the morning + Cups of feed given in the afternoon). \textbf{So based on my findings, the answer is Calculator[(20 * 3) - (6887483 + 25)]}.\\
Act-6: Calculator[(20 * 3) - (6887483 + 25)] \\
Ans: \textbf{-6887448} \\
% \bottomrule
\midrule
\textbf{Ground-Truth Anaswer}: \textbf{-6887448} \\
\bottomrule
\caption{\textsc{False Positive} Examples of generated solutions by IMP-TIP on GSM8K-Hard dataset.  } 
\label{tab:false_pos}
\end{tabular}

\vspace{-10mm}
\end{table*}